# WEAK STATIONARITY OF A MATRIX VALUED DIFFERENTIAL FORM AT SUPERDENSITY POINTS OF ITS VANISHING SET

SILVANO DELLADIO

ABSTRACT. A property of weak stationarity of a matrix valued differential form at superdensity points of its vanishing set is proved. This result is then applied in the context of the Maurer-Cartan equation.

## 1. INTRODUCTION

The main result of this work (cf. Theorem 3.1) establishes a property of weak stationarity of a matrix valued continuous differential form at the superdensity points of its vanishing set. To make this statement more understandable, we now recall very briefly some definitions and properties (referring the reader to Section 2, for a more complete presentation). Let us consider an $M$-dimensional $C^k$ manifold $\mathcal{M}$ and recall that a matrix valued $C^p$ differential $h$-form on $\mathcal{M}$ is a square matrix whose entries are $C^p$ differential $h$-forms on $\mathcal{M}$. The classical formalism for differential forms, i.e., wedge product, exterior differentiation, integration and pullback, extends naturally to matrix valued differential forms (cf. Section 2.2). In this extended formalism it is easy to introduce a notion of distributional exterior derivative, which will be denoted by $\delta$ (cf. Definition 3.1). We also recall that, if $\mathcal{E}$ is a subset of $\mathcal{M}$, then $P \in \mathcal{M}$ is said to be an $m$-density point of $\mathcal{E}$ relative to $\mathcal{M}$ if there is a $C^1$ chart $(\mathcal{W}, \Phi)$ such that $P \in \mathcal{W}$ and

$$\mathcal{L}^M(B_r(\Phi(P)) \setminus \Phi(\mathcal{E} \cap \mathcal{W})) = o(r^m) \qquad (\text{as } r \to 0+),$$

where $\mathcal{L}^M$ and $B_r(\Phi(P))$ are, respectively, the Lebesgue measure on $\mathbb{R}^M$ and the ball of radius $r$ centered at $\Phi(P)$. We observe that this definition does not depend on the choice of the coordinate chart (cf. Section 2.4).

We are now able to state more precisely than before the result in Theorem 3.1: *Let $\mathcal{M}$ be an $M$-dimensional $C^2$ manifold and let $\gamma$ be a matrix valued $C^0$ differential form on $\mathcal{M}$ which has the distributional exterior derivative $\delta\gamma$ of class $C^0$. Then we have $(\delta\gamma)_Q = 0$, whenever $Q$ is an $(M+1)$-density point of $\{P \in \mathcal{M} \,|\, \gamma_P = 0\}$.*

In Section 4, by a simple application of Theorem 3.1, we provide a new proof of the following property in the context of Frobenius theorem about distributions (cf. [5, Theorem 1.3] and [6, Corollary 5.1]): *Let $\mathcal{D}$ be a non-involutive $C^1$ distribution of rank $M$ on a $C^2$ manifold $\mathcal{N}$. Then, for every $M$-dimensional $C^1$ open submanifold $\mathcal{M}$ of $\mathcal{N}$, the tangency set of $\mathcal{M}$ with respect to $\mathcal{D}$ has no $(M+1)$-density points relative to $\mathcal{M}$.*

Section 5 presents an application of Theorem 3.1 in the context of Maurer-Cartan equation. To explain what we are talking about, let us first consider a matrix Lie subgroup $G$ of $\mathrm{Gl}(L, \mathbb{R})$ with Lie algebra $\mathfrak{g}$ and denote its Maurer-Cartan form by $\Gamma_G$. Recall that $\Gamma_G$ is a left-invariant $\mathfrak{g}$-valued smooth differential 1-form on $G$ and
$$d\Gamma_G = -\Gamma_G \wedge \Gamma_G.$$
We have the following well-known theorem, due to Cartan (cf [9, Theorem 1.6.10]): *Let $\mathcal{M}$ be a smooth manifold and let $\phi$ be a $\mathfrak{g}$-valued smooth differential 1-form on $\mathcal{M}$ verifying the Maurer-Cartan equation*
$$(1.1) \qquad d\phi = -\phi \wedge \phi.$$
*Then for all $P \in \mathcal{M}$ there exist a neighborhood $\mathcal{U}$ of $P$ and a smooth map $f : \mathcal{U} \to G$ such that $f^*\Gamma_G = \phi|_{\mathcal{U}}$.*

Relatively to this context, we will provide a structure result for the sets
$$\{P \in \mathcal{U} \mid (f^*\Gamma_G)_P = \phi_P\}$$
under the assumption that $\phi$ does not verify the Maurer-Cartan equation (1.1). In particular, let $\mathcal{M}$ be an $M$-dimensional $C^2$ manifold and let $\phi$ be a $\mathbb{R}^{L \times L}$-valued $C^1$ differential 1-form on $\mathcal{M}$ such that $(d\phi)_Q \neq -(\phi \wedge \phi)_Q$ for all $Q \in \mathcal{M}$. Obviously this condition prevents the possibility of $\phi$ being locally a $C^1$ pullback of $\Gamma_G$ (cf. Remark 5.1). More interesting information on the content of $\{f^*\Gamma_G = \phi|_{\mathcal{U}}\}$ is given in Corollary 5.2, namely: *If $\mathcal{U} \subset \mathcal{M}$ is open and $f : \mathcal{U} \to G$ is a map of class $C^1$, then $\mathcal{U}$ does not contain $(M+1)$-density points of $\{f^*\Gamma_G = \phi|_{\mathcal{U}}\}$.*

## 2. Basic notation and notions

**2.1. Basic notation.** The coordinates of $\mathbb{R}^M$ are denoted by $(x_1, \ldots, x_M)$ so that $dx_1, \ldots, dx_M$ is the standard basis of the dual space of $\mathbb{R}^M$. For simplicity, we set $D_i := \partial/\partial x_i$ and $dx := dx_1 \wedge \cdots \wedge dx_M$. If $p$ is any positive integer not exceeding $M$, then $I(M, p)$ is the family of integer multi-indices $\alpha = (\alpha_1, \ldots, \alpha_p)$ such that $1 \leq \alpha_1 < \cdots < \alpha_p \leq M$. Given a generic map $\Phi : A \to \mathbb{R}^n$ and $v \in \mathbb{R}^n$, we set for simplicity $\{\Phi = v\} := \{P \in A \mid \Phi(P) = v\}$. Let $\mathcal{L}^M$ and $\mathcal{H}^s$ denote, respectively, the Lebesgue measure and the $s$-dimensional Hausdorff measure on $\mathbb{R}^M$. The open ball of radius $r$ centered at $x \in \mathbb{R}^M$ will be denoted by $B_r(x)$. Let $\mathbb{R}^{L \times L}$ be the vector space of all $L \times L$ real matrices and $\mathrm{Gl}(L, \mathbb{R})$ be the Lie group of nondegenerate matrices in $\mathbb{R}^{L \times L}$. The Lie algebra of $\mathrm{Gl}(L, \mathbb{R})$ will be denoted by $\mathfrak{gl}(L, \mathbb{R})$. Since $\mathbb{R}^{L \times L} \simeq \mathbb{R}^{L^2}$ we can denote the natural coordinates on $\mathrm{Gl}(L, \mathbb{R})$ by the matrix notation $(z_{ij})$.

2.2. **Manifolds, differential forms.** In relation to this topic, we will adopt the notations commonly used in the main bibliographic references (see, e.g., [10, 12]). We report here, quickly, just a few of them.

Let $\mathcal{M}$ be an $M$-dimensional $C^k$ manifold. Then a $C^k$ differential $h$-form (resp. $C_c^k$ differential $h$-form, i.e., $C^k$ differential $h$-form with compact support) on $\mathcal{M}$ is a map $\omega : \mathcal{M} \to \Lambda^h T^*\mathcal{M}$ with the following property: If

$$\sum_{\alpha \in I(M,h)} f_\alpha dx_\alpha \qquad (dx_\alpha := dx_{\alpha_1} \wedge \cdots \wedge dx_{\alpha_h})$$

is any local representation of $\omega$, then $f_\alpha$ is of class $C^k$ (resp. $C_c^k$, i.e., $C^k$ with compact support). For any given $P \in \mathcal{M}$, we will use the standard notation $\omega_P$ instead of $\omega(P)$. As we did for real-valued maps, let us set $\{\omega = 0\} := \{P \in \mathcal{M} \,|\, \omega_P = 0\}$ for simplicity. The set of all $C^k$ differential $h$-forms (resp. $C_c^k$ differential $h$-forms) on $\mathcal{M}$ is denoted by $C^k\mathcal{F}^h(\mathcal{M})$ (resp. $C_c^k\mathcal{F}^h(\mathcal{M})$).

Let $\mathcal{M}$ be a $C^k$ imbedded submanifold of a $C^k$ manifold $\mathcal{N}$ and let $\iota : \mathcal{M} \hookrightarrow \mathcal{N}$ be the inclusion map. If $\omega \in C^{k-1}\mathcal{F}^h(\mathcal{N})$, then the $C^{k-1}$ differential $h$-form $\iota^*\omega$ (i.e., the restriction of $\omega$ to $\mathcal{M}$) will be denoted by $\omega|_{\mathcal{M}}$.

We also need matrix-valued differential forms, i.e., matrices whose entries are differential forms. If $\mathcal{M}$ is a $C^k$ manifold and $L$ is a positive integer then $\mathrm{Mat}_L C^p\mathcal{F}^h(\mathcal{M})$ is the set of all $L \times L$ matrices

$$(\omega^{(ij)}) = \begin{pmatrix} \omega^{(11)} & \cdots & \omega^{(1L)} \\ \vdots & \ddots & \vdots \\ \omega^{(L1)} & \cdots & \omega^{(LL)} \end{pmatrix}, \text{ with } \omega^{(ij)} \in C^p\mathcal{F}^h(\mathcal{M}).$$

For the sake of convenience, we will sometimes (e.g. in Section 5 below) refer to the members of $\mathrm{Mat}_L C^p\mathcal{F}^h(\mathcal{M})$ by simply calling them $C^p$ differential $h$-forms as well. The subset of $\mathrm{Mat}_L C^p\mathcal{F}^h(\mathcal{M})$ whose members have all the entries in $C_c^p\mathcal{F}^h(\mathcal{M})$ is denoted by $\mathrm{Mat}_L C_c^p\mathcal{F}^h(\mathcal{M})$. If $\omega = (\omega^{(ij)}) \in \mathrm{Mat}_L C_c^p\mathcal{F}^h(\mathcal{M})$ then we set $\mathrm{supp}(\omega) := \cup_{i,j} \mathrm{supp}(\omega_{ij})$.

If $\omega = (\omega^{(ij)}) \in \mathrm{Mat}_L C^p\mathcal{F}^h(\mathcal{M})$, then we define

$$\omega_P := (\omega_P^{(ij)}), \quad \omega_P(v_1, \ldots, v_h) := (\omega_P^{(ij)}(v_1, \ldots, v_h))$$

for all $P \in \mathcal{M}$ and $v_1, \ldots, v_h \in T_P\mathcal{M}$. If $p \geq 1$, we define the exterior differentiation $d : \mathrm{Mat}_L C^p\mathcal{F}^h(\mathcal{M}) \to \mathrm{Mat}_L C^{p-1}\mathcal{F}^{h+1}(\mathcal{M})$ by

$$d(\omega^{(ij)}) := (d\omega^{(ij)}).$$

Observe that $d$ is linear and $d \circ d = 0$. If $\mathcal{N}$ is another $C^k$ manifold and $f : \mathcal{M} \to \mathcal{N}$ is a $C^p$ map, the pullback

$$f^* : \mathrm{Mat}_L C^p\mathcal{F}^h(\mathcal{N}) \to \mathrm{Mat}_L C^{p-1}\mathcal{F}^h(\mathcal{M})$$

is defined as follows

$$f^*(\omega^{(ij)}) := (f^*\omega^{(ij)}).$$

The exterior product of of two matrix-valued differential forms
$$\lambda = (\lambda^{(ij)}) \in \mathrm{Mat}_L C^p \mathcal{F}^l(\mathcal{M}), \quad \mu = (\mu^{(ij)}) \in \mathrm{Mat}_L C^p \mathcal{F}^m(\mathcal{M})$$
is the matrix-valued differential form $\lambda \wedge \mu \in \mathrm{Mat}_L C^p \mathcal{F}^{l+m}(\mathcal{M})$ whose entries are defined by
$$(\lambda \wedge \mu)^{(ij)} := \sum_{q=1}^{L} \lambda^{(iq)} \wedge \mu^{(qj)}.$$
A trivial computation shows that differentiating the exterior product of matrix-valued differential forms yields the usual formula (provided $k \geq 1$):
$$d(\lambda \wedge \mu) = d\lambda \wedge \mu + (-1)^l \lambda \wedge d\mu.$$
A matrix-valued differential form $\omega = (\omega^{(ij)}) \in \mathrm{Mat}_L C^0 \mathcal{F}^M(\mathcal{M})$ is said to be integrable on $\mathcal{M}$ if every $\omega^{(ij)}$ is integrable on $\mathcal{M}$. In this case we set
$$\int_{\mathcal{M}} \omega := \left( \int_{\mathcal{M}} \omega^{(ij)} \right).$$

Let us recall that a $C^1$ Riemannian manifold $(\mathcal{M}, g)$ with the associated Riemannian distance function is a metric space whose topology coincides to the original manifold topology, cf. [10, Theorem 13.29]. Hence one can define the corresponding $s$-dimensional Hausdorff measure $\mathcal{H}_g^s$, cf. [8, Section 2.10.2], [13, Chapter 12]. The open metric ball of radius $r$ centered at $P \in \mathcal{M}$ will be denoted by $\mathcal{B}_g(P, r)$.

2.3. **Hausdorff measure on manifolds.** For the convenience of the reader, we recall the following well-known properties of the Hausdorff measure $\mathcal{H}_g^s$ on a $C^1$ Riemannian manifold $(\mathcal{N}, g)$:

- If $s = \dim \mathcal{N}$, then $\mathcal{H}_g^s(B) = V_g(B)$ for all Borel sets $B \subset \mathcal{N}$, where $V_g$ denotes the standard volume form of $(\mathcal{N}, g)$, cf. [8, Section 3.2.46], [13, Proposition 12.6].
- If $\mathcal{M}$ is a $C^1$ imbedded submanifold of $\mathcal{N}$ and $g_{\mathcal{M}}$ denotes the induced metric, then one has $\mathcal{H}_{g_{\mathcal{M}}}^s(B) = \mathcal{H}_g^s(B)$ for all Borel sets $B \subset \mathcal{M}$, cf. [13, Proposition 12.7].
- If $g$ denotes the standard Euclidean metric on $\mathbb{R}^N$, then one obviously has $\mathcal{H}_g^s = \mathcal{H}^s$. In particular, $\mathcal{H}_g^N$ is the $N$-dimensional Lebesgue measure.

Another property which follows readily from [8, Section 3.2.46] is this one.

**Proposition 2.1.** *Let $\mathcal{N}$ be a $C^1$ manifold, $\mathcal{E} \subset \mathcal{N}$ and $s \in [0, +\infty)$. The following are equivalent:*

(1) *For every $C^1$ chart $(\mathcal{W}, \Phi)$ of $\mathcal{N}$, one has $\mathcal{H}^s(\Phi(\mathcal{W} \cap \mathcal{E})) = 0$.*
(2) *For every $C^1$ Riemannian metric $g$ on $\mathcal{N}$, one has $\mathcal{H}_g^s(\mathcal{E}) = 0$.*

  (3) *There exists a $C^1$ Riemannian metric $g$ on $\mathcal{N}$ such that $\mathcal{H}_g^s(\mathcal{E}) = 0$.*

2.4. **Superdensity.** Also the following proposition is a consequence of [8, Section 3.2.46], cf. [5, Proposition 3.3].

**Proposition 2.2.** *Let $\mathcal{N}$ be a $N$-dimensional $C^1$ manifold, $\mathcal{E} \subset \mathcal{N}$, $P \in \mathcal{N}$ and $m \in [N, +\infty)$. The following are equivalent:*

  (1) *There is a $C^1$ chart $(\mathcal{W}, \Phi)$ of $\mathcal{N}$ such that $P \in \mathcal{W}$ and*
  $$\mathcal{L}^N(B_r(\Phi(P)) \setminus \Phi(\mathcal{E} \cap \mathcal{W})) = o(r^m) \qquad (as\ r \to 0+).$$

  (2) *For every $C^1$ Riemannian metric $g$ on $\mathcal{N}$, one has*
  $$\mathcal{H}_g^N(\mathcal{B}_g(P, r) \setminus \mathcal{E}) = o(r^m) \qquad (as\ r \to 0+).$$

  (3) *There exists a $C^1$ Riemannian metric $g$ on $\mathcal{N}$ such that*
  $$\mathcal{H}_g^N(\mathcal{B}_g(P, r) \setminus \mathcal{E}) = o(r^m) \qquad (as\ r \to 0+).$$

**Definition 2.1.** *If any or, equivalently, all of the conditions of Proposition 2.2 are satisfied, then we say that $P$ is an $m$-density point of $\mathcal{E}$ (relative to $\mathcal{N}$). The set of all $m$-density points of $\mathcal{E}$ is denoted by $\mathcal{E}^{(m)}$, cf. [5].*

**Remark 2.1.** *Let $\mathcal{N}$ and $\mathcal{E}$ be as in Proposition 2.2. The following facts occur:*

- *Every interior point of $\mathcal{E}$ is an $m$-density point of $\mathcal{E}$, for all $m \in [N, +\infty)$. Thus, whenever $\mathcal{E}$ is open, one has $\mathcal{E} \subset \mathcal{E}^{(m)}$ for all $m \in [N, +\infty)$.*
- *If $N \leq m_1 \leq m_2 < +\infty$, then $\mathcal{E}^{(m_2)} \subset \mathcal{E}^{(m_1)}$. In particular, one has $\mathcal{E}^{(m)} \subset \mathcal{E}^{(N)}$ for all $m \in [N, +\infty)$.*
- *Let $\{\mathcal{E}_j\}_{j \in J}$ be any family of subsets of $\mathcal{N}$ and $m \in [N, +\infty)$.*
  - *One has*
  $$\left(\bigcap_{j \in J} \mathcal{E}_j\right)^{(m)} \subset \bigcap_{j \in J} \mathcal{E}_j^{(m)};$$
  - *If $J$ is finite, then*

(2.1)
  $$\left(\bigcap_{j \in J} \mathcal{E}_j\right)^{(m)} = \bigcap_{j \in J} \mathcal{E}_j^{(m)};$$

  - *If $J$ is countable infinite, then (2.1) can fail to be true, e.g., $\mathcal{N} = \mathbb{R}^2$ and*
  $$\mathcal{E}_j := B_{1/j}(O) \qquad (j = 1, 2, \dots).$$

**Remark 2.2.** *For convenience of the reader, we recall some known results in the special case when $\mathcal{N} = \mathbb{R}^N$ (which actually could be easily generalized):*

- *If $E \subset \mathbb{R}^N$ is $\mathcal{L}^N$-measurable then: $x \in E^{(N)}$ if and only if $x$ is a Lebesgue density point of $E$, hence $\mathcal{L}^N(E \Delta E^{(N)}) = 0$. In particular, it follows that $(E^{(N)})^{(N)} = E^{(N)}$.*

- If $E \subset \mathbb{R}^N$, then $E^{(m)}$ is $\mathcal{L}^N$-measurable, for all $m \in [N, +\infty)$ (cf. [3, Proposition 3.1]).
- Every open set $U \subset \mathbb{R}^N$ can be approximated in measure by uniformly $N$-dense closed subsets of $\overline{U}$. More precisely: For all $C < \mathcal{L}^N(U)$ there exists a closed set $F \subset \overline{U}$ such that $\mathcal{L}^N(F) > C$ and $F^{(m)} = \emptyset$ for all $m > N$ (obviously one has $F^{(N)} \subset F$ and $\mathcal{L}^N(F \setminus F^{(N)}) = 0$), cf. [4, Proposition 5.4].
- Let $N \geq 2$ and $E \subset \mathbb{R}^N$ be a set of finite perimeter, so that $\mathcal{H}^{N-1}(\partial^* E) < +\infty$ (where $\partial^* E$ is the reduced boundary of $E$, cf. [11, Theorem 15.9]). Then $\mathcal{L}^N(E \setminus E^{(m_0)}) = 0$, with

$$m_0 := N + 1 + \frac{1}{N-1},$$

cf. Theorem 1 in [7, Section 6.1.1] (compare also [2, Lemma 4.1]). Moreover, the number $m_0$ is the maximum order of density common to all sets of finite perimeter. More precisely, the following property holds (cf. [3, Proposition 4.1]): For all $m > m_0$ there exists a compact set $F_m$ of finite perimeter in $\mathbb{R}^N$ such that $\mathcal{L}^N(F_m) > 0$ and $F_m^{(m)} = \emptyset$.

## 3. The main result

Throughout this section $\mathcal{M}$ and $k$ will denote, respectively, a $M$-dimensional manifold and the regularity class of $\mathcal{M}$. We will assume $k \geq 1$, if not otherwise stated.

**Remark 3.1.** *Let $l \leq M$ and $\lambda \in Mat_L C^0 \mathcal{F}^l(\mathcal{M})$. Then $\lambda = 0$ if and only if*

$$\int_{\mathcal{M}} \lambda \wedge \mu = 0$$

*for all $\mu \in Mat_L C_c^k \mathcal{F}^{M-l}(\mathcal{M})$.*

From Remark 3.1 we get immediately the following proposition.

**Proposition 3.1.** *Let $\lambda \in Mat_L C^0 \mathcal{F}^h(\mathcal{M})$, with $h \leq M-1$, satisfy the following property: there exists $\mu \in Mat_L C^0 \mathcal{F}^{h+1}(\mathcal{M})$ such that $\int_{\mathcal{M}} \lambda \wedge d\varphi = \int_{\mathcal{M}} \mu \wedge \varphi$, for all $\varphi \in Mat_L C_c^k \mathcal{F}^{M-h-1}(\mathcal{M})$. Then $\mu$ is uniquely determined.*

**Definition 3.1.** *Let the assumptions of Proposition 3.1 be verified. Then we say that $\lambda$ has the distributional exterior derivative (DED) in $Mat_L C^0 \mathcal{F}^{h+1}(\mathcal{M})$. The latter is defined as $\delta\lambda := (-1)^{h+1}\mu$, so that*

$$\tag{3.1} \int_{\mathcal{M}} \lambda \wedge d\varphi = (-1)^{h+1} \int_{\mathcal{M}} \delta\lambda \wedge \varphi$$

*for all $\varphi \in Mat_L C_c^k \mathcal{F}^{M-h-1}(\mathcal{M})$.*

**Remark 3.2.** *Let $l$ be an integer such that $1 \leq l \leq k$. Then a standard approximation argument shows that $Mat_L C_c^k \mathcal{F}^{M-h-1}(\mathcal{M})$ is dense in $Mat_L C_c^l \mathcal{F}^{M-h-1}(\mathcal{M})$, with respect to the $C^l$ topology. Hence in Definition 3.1 we can equivalently assume that (3.1) holds for all $\varphi \in Mat_L C_c^l \mathcal{F}^{M-h-1}(\mathcal{M})$.*

The following propositions state some expected properties. We observe that the first three are trivial.

**Proposition 3.2.** *If $\lambda \in Mat_L C^0 \mathcal{F}^h(\mathcal{M})$ has the DED in $Mat_L C^0 \mathcal{F}^{h+1}(\mathcal{M})$ and $U \subset \mathcal{M}$ is open, then $\lambda|_U$ has the DED in $Mat_L C^0 \mathcal{F}^{h+1}(U)$ and $\delta(\lambda|_U) = (\delta\lambda)|_U$.*

**Proposition 3.3.** *If $\lambda \in Mat_L C^1 \mathcal{F}^h(\mathcal{M})$ then $\lambda$ has the DED in $Mat_L C^0 \mathcal{F}^{h+1}(\mathcal{M})$ and $\delta\lambda = d\lambda$.*

**Proposition 3.4.** *Let $\lambda, \mu \in Mat_L C^0 \mathcal{F}^h(\mathcal{M})$ have the DED in $Mat_L C^0 \mathcal{F}^{h+1}(\mathcal{M})$. Then, for all $a, b \in \mathbb{R}$, the matrix-valued differential form $a\lambda + b\mu \in Mat_L C^0 \mathcal{F}^h(\mathcal{M})$ has the DED in $Mat_L C^0 \mathcal{F}^{h+1}(\mathcal{M})$ and $\delta(a\lambda + b\mu) = a\,\delta\lambda + b\,\delta\mu$.*

**Proposition 3.5.** *Let $\mathcal{M}$ be of class $C^k$, with $k \geq 2$. If $\lambda \in Mat_L C^0 \mathcal{F}^h(\mathcal{M})$ has the DED in $Mat_L C^0 \mathcal{F}^{h+1}(\mathcal{M})$, then $\delta\lambda$ has the DED in $Mat_L C^0 \mathcal{F}^{h+2}(\mathcal{M})$ and $\delta(\delta\lambda) = 0$.*

*Proof.* Let $\lambda \in \mathrm{Mat}_L C^0 \mathcal{F}^h(\mathcal{M})$ have the DED in $\mathrm{Mat}_L C^0 \mathcal{F}^{h+1}(\mathcal{M})$. Then, by Definition 3.1 and Remark 3.2 (with $l = k - 1$), we obtain
$$\int_{\mathcal{M}} \delta\lambda \wedge d\varphi = (-1)^{h+1} \int_{\mathcal{M}} \lambda \wedge d(d\varphi) = 0 = (-1)^{h+2} \int_{\mathcal{M}} 0 \wedge \varphi$$
for all $\varphi \in \mathrm{Mat}_L C^k \mathcal{F}^{M-h-2}(\mathcal{M})$. □

**Remark 3.3.** *Combining Proposition 3.3 and Proposition 3.5, we obtain the following property: If $k \geq 2$ and $\lambda \in Mat_L C^1 \mathcal{F}^h(\mathcal{M})$, then $d\lambda$ has the DED in $Mat_L C^0 \mathcal{F}^{h+2}(\mathcal{M})$ and $\delta(d\lambda) = 0$.*

**Proposition 3.6.** *Let $\mathcal{M}$ be of class $C^k$, with $k \geq 2$. Moreover consider a $C^2$ manifold $\mathcal{N}$, a $C^1$ map $f : \mathcal{M} \to \mathcal{N}$ and $\omega \in Mat_L C^1 \mathcal{F}^h(\mathcal{N})$, with $h \leq M - 1$. Then $f^*\omega$ has the DED in $Mat_L C^0 \mathcal{F}^{h+1}(\mathcal{M})$ and $\delta(f^*\omega) = f^*(d\omega)$.*

*Proof.* Consider $\varphi \in \mathrm{Mat}_L C_c^k \mathcal{F}^{M-h-1}(\mathcal{M})$. Then for all $x \in \mathrm{supp}(\varphi)$ there exists an open set $\mathcal{V}^{(x)} \subset \mathcal{M}$ and a countable family $\{f_j^{(x)}\} \subset C^2(\mathcal{V}^{(x)}, \mathcal{N})$ such that $f_j^{(x)} \to f$ (as $j \to \infty$) with respect to the $C^1(\mathcal{V}^{(x)}, \mathcal{N})$ topology. Since $\mathrm{supp}(\varphi)$ is compact, there exists a finite set $\{x_1, \ldots, x_N\} \subset \mathrm{supp}(\varphi)$ such that
$$\mathrm{supp}(\varphi) \subset \mathcal{V} := \cup_i \mathcal{V}^{(x_i)}. \tag{3.2}$$
By [12, Theorem 2.2.14] we can find $\{\eta_1, \ldots, \eta_N\} \subset C^2(\mathcal{M})$ such that
$$\eta_i \geq 0, \quad \mathrm{supp}(\eta_i) \subset \mathcal{V}^{(x_i)}, \quad \sum_i \eta_i|_{\mathcal{V}} = 1.$$
If we extend every $f_j^{(x_i)}$ arbitrarily to all of $\mathcal{M}$ and define
$$f_j := \sum_i \eta_i f_j^{(x_i)} \in C_c^2(\mathcal{M}, \mathcal{N}) \qquad (j = 1, 2, \ldots)$$
then $f_j|_{\mathcal{V}} \to f|_{\mathcal{V}}$ (as $j \to \infty$) with respect to the $C^1(\mathcal{V}, \mathcal{N})$ topology. Moreover we have
$$\int_{\mathcal{V}} f_j^*(d\omega) \wedge \varphi = \int_{\mathcal{V}} d(f_j^*\omega) \wedge \varphi = (-1)^{h+1} \int_{\mathcal{V}} (f_j^*\omega) \wedge d\varphi.$$

Hence, letting $j \to +\infty$, we obtain

$$\int_{\mathcal{V}} f^*(d\omega) \wedge \varphi = (-1)^{h+1} \int_{\mathcal{V}} (f^*\omega) \wedge d\varphi$$

that is (by (3.2))

$$\int_{\mathcal{M}} f^*(d\omega) \wedge \varphi = (-1)^{h+1} \int_{\mathcal{M}} (f^*\omega) \wedge d\varphi.$$

The conclusion follows from the arbitrariness of $\varphi$. □

Let us now state and prove the main result.

**Theorem 3.1.** *Let $\mathcal{M}$ be of class $C^k$, with $k \geq 2$. Moreover let $h \leq M - 1$ and consider $\gamma \in \mathrm{Mat}_L C^0 \mathcal{F}^h(\mathcal{M})$ which has the DED in $\mathrm{Mat}_L C^0 \mathcal{F}^{h+1}(\mathcal{M})$. If define*

$$\mathcal{Z}_\gamma := \{P \in \mathcal{M} \,|\, \gamma_P = 0\}$$

*then $(\delta\gamma)_Q = 0$ for all $Q \in \mathcal{Z}_\gamma^{(M+1)}$.*

*Proof.* First of all, set for simplicity $B_r := B_r(0) \subset \mathbb{R}^M$ and let $\rho \in (0,1)$. Then consider $g \in C_c^2(B_1)$ such that $0 \leq g \leq 1$, $g|_{B_\rho} \equiv 1$ and

$$|D_i g| \leq \frac{2}{1-\rho} \qquad (i = 1, \ldots, M).$$

For $r > 0$, define $g_r \in C_c^2(B_r)$ as

$$g_r(x) := g\left(\frac{x}{r}\right), \quad x \in B_r$$

and observe that (for all $x \in B_r$ and $i = 1, \ldots, M$)

$$(3.3) \qquad |D_i g_r(x)| = \frac{1}{r}\left|D_i g\left(\frac{x}{r}\right)\right| \leq \frac{2}{r(1-\rho)}.$$

Now consider an arbitrary $Q \in \mathcal{Z}_\gamma^{(M+1)}$ and let $(\mathcal{U}, \Phi)$ be a $C^2$ coordinate chart on $\mathcal{M}$ such that $Q \in \mathcal{U}$ and $\Phi(Q) = 0 \in \mathbb{R}^M$. Observe that

$$(3.4) \qquad \mathcal{L}^M(B_r \setminus \Phi(Z_\gamma)) = o(r^{M+1}) \qquad (\text{as } r \to 0+)$$

by Definition 2.1.

Now set for simplicity $U := \Phi(\mathcal{U})$ and let $\theta \in \mathrm{Mat}_L C^2 \mathcal{F}^{M-1-h}(U)$ be chosen arbitrarily. Obviously there must be $(F_\theta^{(ij)}) \in \mathrm{Mat}_L C^0 \mathcal{F}^0(U)$ such that

$$(3.5) \qquad [(\Phi^{-1})^*(\delta\gamma)] \wedge \theta = (F_\theta^{(ij)} \, dx),$$

hence, for all $i, j$, we have (provided $r$ is small enough)

$$\left| \int_{B_r} g_r \, F_\theta^{(ij)} \, dx \right| = \left| \int_{B_r} g_r \left( [(\Phi^{-1})^*(\delta\gamma)] \wedge \theta \right)^{(ij)} \right|$$

$$= \left| \int_{\Phi^{-1}(B_r)} (g_r \circ \Phi) \left( (\delta\gamma) \wedge (\Phi^*\theta) \right)^{(ij)} \right|$$

$$= \left| \int_{\Phi^{-1}(B_r)} \left( (\delta\gamma) \wedge [(g_r \circ \Phi) \, \Phi^*\theta] \right)^{(ij)} \right|$$

$$= \left| \int_{\Phi^{-1}(B_r)} \left( \gamma \wedge d[(g_r \circ \Phi) \, \Phi^*\theta] \right)^{(ij)} \right|$$

$$\leq \left| \int_{\Phi^{-1}(B_r) \setminus \mathcal{Z}_\gamma} \left( \gamma \wedge d(g_r \circ \Phi) \wedge \Phi^*\theta \right)^{(ij)} \right|$$

$$+ \left| \int_{\Phi^{-1}(B_r) \setminus \mathcal{Z}_\gamma} (g_r \circ \Phi) \left( \gamma \wedge \Phi^*(d\theta) \right)^{(ij)} \right|$$

$$= \left| \int_{B_r \setminus \Phi(\mathcal{Z}_\gamma)} \left( [(\Phi^{-1})^*\gamma] \wedge dg_r \wedge \theta \right)^{(ij)} \right|$$

$$+ \left| \int_{B_r \setminus \Phi(\mathcal{Z}_\gamma)} g_r \left( [(\Phi^{-1})^*\gamma] \wedge d\theta \right)^{(ij)} \right|.$$

Recalling (3.3), we obtain

$$\left| \int_{B_r} g_r \, F_\theta^{(ij)} \, dx \right| \leq C \, \mathcal{L}^M(B_r \setminus \Phi(\mathcal{Z}_\gamma)) \left( \frac{1}{r(1-\rho)} + 1 \right).$$

On the other hand, the triangle inequality yields

$$\left| \int_{B_r} g_r \, F_\theta^{(ij)} \, dx \right| \geq \left| \int_{B_{\rho r}} g_r \, F_\theta^{(ij)} \, dx \right| - \left| \int_{B_r \setminus B_{\rho r}} g_r \, F_\theta^{(ij)} \, dx \right|$$

$$= \left| \int_{B_{\rho r}} F_\theta^{(ij)} \, dx \right| - \left| \int_{B_r \setminus B_{\rho r}} g_r \, F_\theta^{(ij)} \, dx \right|.$$

It follows that

$$\rho^M \left| \fint_{B_{\rho r}} F_\theta^{(ij)} \, dx \right| \leq \frac{C \mathcal{L}^M(B_r \setminus \Phi(\mathcal{Z}_\gamma))}{r^M} \left( \frac{1}{r(1-\rho)} + 1 \right) + \frac{C(r^M - \rho^M r^M)}{r^M}$$

$$= \frac{C \mathcal{L}^M(B_r \setminus \Phi(\mathcal{Z}_\gamma))}{r^{M+1}} \left( \frac{1}{1-\rho} + r \right) + C(1 - \rho^M).$$

Then, by first letting $r \to 0+$ (and recalling (3.4)) and then letting $\rho \to 1-$, we obtain $F_\theta^{(ij)}(0) = 0$ (for all $i, j$). The conclusion follows from the identity (3.5) and the arbitrariness of $\theta$. □

The following simple corollary of Theorem 3.1 will be useful below.

**Corollary 3.1.** *Let $\mathcal{M}$ and $\mathcal{N}$ be two $C^2$ manifolds, let $f : \mathcal{M} \to \mathcal{N}$ be a $C^1$ map and $\omega \in Mat_L C^1 \mathcal{F}^h(\mathcal{N})$, with $h + 1 \leq M := \dim \mathcal{M}$. Moreover consider $\mu \in Mat_L C^0 \mathcal{F}^h(\mathcal{M})$ which has the DED in $Mat_L C^0 \mathcal{F}^{h+1}(\mathcal{M})$ and define*

$$\mathcal{A}_{f,\omega,\mu} := \{P \in \mathcal{M} \mid \mu_P = (f^*\omega)_P\}.$$

*Then $(\delta\mu)_Q = (f^* d\omega)_Q$, for all $Q \in \mathcal{A}_{f,\omega,\mu}^{(M+1)}$.*

*Proof.* Define $\gamma := \mu - f^*\omega \in \mathrm{Mat}_L C^0 \mathcal{F}^h(\mathcal{M})$ and observe that $\mathcal{A}_{f,\omega,\mu} = \mathcal{Z}_\gamma$, hence

$$\mathcal{A}_{f,\omega,\mu}^{(M+1)} = \mathcal{Z}_\gamma^{(M+1)}.$$

Moreover, by Proposition 3.4 and Proposition 3.6, the form $\gamma$ has the distributional exterior derivative in $\mathrm{Mat}_L C^0 \mathcal{F}^{h+1}(\mathcal{M})$ and

$$\delta\gamma = \delta\mu - f^* d\omega.$$

The conclusion follows from Theorem 3.1. □

## 4. Applications I

From Corollary 3.1 we can easily derive [6, Theorem 3.1], which states a low-density property for the integral set of a submanifold with respect to a non-integrable exterior differential system. Before showing this application, let us briefly set the context. Consider a $C^2$ manifold $\mathcal{N}$ and an arbitrary family $\mathcal{O}$ of $C^1$ differential forms on $\mathcal{N}$. Moreover let $f : U \subset \mathbb{R}^M \to \mathcal{N}$ (where $U$ is open), be any imbedding of class $C^1$ and define

$$\mathcal{I}(f, \mathcal{O}) := \bigcap_{\omega \in \mathcal{O}} \{f^*\omega = 0\}.$$

Then [6, Theorem 3.1] states that

$$U \cap \mathcal{I}(f, \mathcal{O})^{(M+1)} \subset \bigcap_{\omega \in \mathcal{O}} \{f^* d\omega = 0\}.$$

Now let $V_M(\mathcal{O})_y$ denote the set of all $M$-dimensional integral elements of $\mathcal{O}$ at $y \in \mathcal{N}$ (cf. Definition 1.1 in Section 1 of [1, Chapter III] and the first definition in Section 1 of [14, Chapter III]) and assume that

(4.1) *For all $y \in \mathcal{N}$ and $\Sigma \in V_M(\mathcal{O})_y$ there is $\omega \in \mathcal{O}$ such that $(d\omega)_y|_\Sigma \neq 0$.*

We naturally expect that condition (4.1) prevents the existence of interior points in $\mathcal{I}(f, \mathcal{O})$, but the structure of $\mathcal{I}(f, \mathcal{O})$ can be described more precisely by using the notion of superdensity. Indeed in [6, Corollary 3.2], which follows trivially from [6, Theorem 3.1], we have proved that one has

(4.2) $$U \cap \mathcal{I}(f, \mathcal{O})^{(M+1)} = \emptyset.$$

We can finally apply Corollary 3.1 to prove the following result, which in turn served to prove [6, Theorem 3.1] very easily.

**Theorem 4.1** (Theorem 3.2 of [6])**.** *Let $\omega \in C^1\mathcal{F}^h(\mathcal{N})$ and $f : U \subset \mathbb{R}^M \to \mathcal{N}$ (where $U$ is open) be a $C^1$ map. Then*
$$U \cap \{d\lambda = f^*\omega\}^{(M+1)} \subset \{f^*d\omega = 0\}$$
*for every $\lambda \in C^1\mathcal{F}^{h-1}(U)$.*

*Proof.* Observe that $\mu := (d\lambda) \in \text{Mat}_1 C^0 \mathcal{F}^h(U)$ has the DED in $\text{Mat}_1 C^0 \mathcal{F}^{h+1}(U)$ and $\delta\mu = 0$, by Remark 3.3. Hence and by Corollary 3.1 (with $L = 1$) we get $(f^*d\omega)_Q = 0$ for all $Q \in \mathcal{A}_{f,(\omega),\mu}^{(M+1)} = U \cap \{d\lambda = f^*\omega\}^{(M+1)}$. □

**Remark 4.1.** *If $\mathcal{O}$ is a family of linearly independent $C^1$ differential 1-forms defining a distribution $\mathcal{D}$ of rank $M$ on $\mathcal{N}$ (cf.[10, Chapter 19]), then, for all $y \in \mathcal{N}$, the $M$-plane $\mathcal{D}_y$ is the only $M$-dimensional integral element of $\mathcal{O}$ at $y$, i.e., $V_M(\mathcal{O})_y = \{\mathcal{D}_y\}$. Hence:*

- *The set $\mathcal{I}(f,\mathcal{O})$ coincides with the tangency set of $f(U)$ with respect to $\mathcal{D}$;*
- *The condition (4.1) is verified if and only if $\mathcal{D}$ is non-involutive at each point of $\mathcal{N}$, cf. [10, Proposition 19.8].*

*Thus the structure identity (4.2) proves that if $\mathcal{D}$ is non-involutive at each point of $\mathcal{N}$ then the following property holds: For every $M$-dimensional $C^1$ open submanifold $\mathcal{M}$ of $\mathcal{N}$, the tangency set of $\mathcal{M}$ with respect to $\mathcal{D}$ has no $(M+1)$-density points relative to $\mathcal{M}$, cf. [5, Theorem 1.3] and [6, Corollary 5.1].*

5. APPLICATIONS II, THE CONTEXT OF MAURER-CARTAN FORM

Let us consider any matrix Lie subgroup $G$ of $\text{Gl}(L,\mathbb{R})$ with Lie algebra $\mathfrak{g} \subset \mathfrak{gl}(L,\mathbb{R})$ and let $\iota : G \to \text{Gl}(L,\mathbb{R})$ be the inclusion map. Then let $\gamma \in \text{Mat}_L C^\infty \mathcal{F}^1(\text{Gl}(L,\mathbb{R}))$ be defined at $z = (z_{ij}) \in \text{Gl}(L,\mathbb{R})$ as
$$\gamma_z := (z_{ij})^{-1}(dz_{ij})$$
and define the Maurer-Cartan form of $G$ as
$$\Gamma_G := \iota^*\gamma \in \text{Mat}_L C^\infty \mathcal{F}^1(G).$$
Observe that $\gamma$ is the Maurer-Cartan form of $\text{Gl}(L,\mathbb{R})$. Recall that $\Gamma_G$ is left-invariant, takes values in $\mathfrak{g}$ and satisfies the Maurer-Cartan equation, that is

(5.1) $$d\Gamma_G = -\Gamma_G \wedge \Gamma_G,$$

cf. [9, Section 1.6].

**Remark 5.1.** *Consider a $C^2$ manifold $\mathcal{M}$, $\phi \in \text{Mat}_L C^1 \mathcal{F}^1(\mathcal{M})$ and assume that the following property holds: For all $P \in \mathcal{M}$ there exist a neighborhood $\mathcal{U}$ of $P$ and a $C^1$ map $f : \mathcal{U} \to G$ such that $f^*\Gamma_G = \phi|_\mathcal{U}$. Then, first of all, $\phi$ takes values in $\mathfrak{g}$. Moreover, by Proposition 3.3, Proposition 3.6 and (5.1), one has*
$$d(f^*\Gamma_G) = \delta(f^*\Gamma_G) = f^*(d\Gamma_G) = -f^*(\Gamma_G \wedge \Gamma_G) = -(f^*\Gamma_G) \wedge (f^*\Gamma_G)$$

*that is*
$$(d\phi)|_{\mathcal{U}} = -(\phi \wedge \phi)|_{\mathcal{U}}.$$

*Relative to the opposite implication, it is well known that a $\mathfrak{g}$-valued smooth differential 1-form satisfying the Maurer-Cartan equation is always, at least locally, a smooth pullback of the Maurer-Cartan form. In fact the following theorem holds, cf [9, Theorem 1.6.10].*

**Theorem 5.1** (Cartan). *Let $\mathcal{M}$ be a smooth manifold and let $\phi$ be a $\mathfrak{g}$-valued smooth differential 1-form on $\mathcal{M}$ satisfying the identity $d\phi = -\phi \wedge \phi$. Then for all $P \in \mathcal{M}$ there exist a neighborhood $\mathcal{U}$ of $P$ and a smooth map $f : \mathcal{U} \to G$ such that $f^*\Gamma_G = \phi|_{\mathcal{U}}$. Moreover, if $f_1, f_2 : \mathcal{U} \to G$ are any two smooth maps with this property, then there exists $a \in G$ such that $f_2(Q) = af_1(Q)$ for all $Q \in \mathcal{U}$.*

Remark 5.1 shows that, if $\mathcal{M}$ is a $C^2$ manifold and $\phi \in \text{Mat}_L C^1 \mathcal{F}^1(\mathcal{M})$, the occurrence of condition

(5.2) $$(d\phi)_Q \neq -(\phi \wedge \phi)_Q, \text{ for all } Q \in \mathcal{M}$$

prevents the possibility of $\phi$ being locally a $C^1$ pullback of the Maurer-Cartan form $\Gamma_G$. Thus, whatever the choice of $C^1$ map $f : \mathcal{U} \subset \mathcal{M} \to G$, the set $\{f^*\Gamma_G = \phi|_{\mathcal{U}}\}$ cannot have interior points. In Corollary 5.2 below we provide a structure result for this set, under assumption (5.2), by using superdensity.

Now we provide an application of Corollary 3.1, which is the natural counterpart in this context of Theorem 4.1 in Section 4.

**Theorem 5.2.** *Let $\mathcal{M}$ be a $M$-dimensional $C^2$ manifold and let $\phi \in \text{Mat}_L C^0 \mathcal{F}^1(\mathcal{M})$ have the DED in $\text{Mat}_L C^0 \mathcal{F}^2(\mathcal{M})$. Moreover, let $\mathcal{U} \subset \mathcal{M}$ be open and consider a $C^1$ map $f : \mathcal{U} \to G$. Then $(\delta\phi)_Q = -(\phi \wedge \phi)_Q$ for all $Q \in \mathcal{U} \cap \{f^*\Gamma_G = \phi|_{\mathcal{U}}\}^{(M+1)}$.*

*Proof.* Let $Q \in \mathcal{U} \cap \{f^*\Gamma_G = \phi|_{\mathcal{U}}\}^{(M+1)}$ and observe that

(5.3) $$(f^*\Gamma_G)_Q = \phi_Q,$$

by continuity. We observe also that, by Proposition 3.2, $\phi|_{\mathcal{U}}$ has the DED in $\text{Mat}_L C^0 \mathcal{F}^2(\mathcal{U})$ and $\delta(\phi|_{\mathcal{U}}) = (\delta\phi)|_{\mathcal{U}}$. If we now apply Corollary 3.1 with

$$\mathcal{M} := \mathcal{U}, \quad \mathcal{N} := G, \quad \omega := \Gamma_G, \quad \mu := \phi|_{\mathcal{U}},$$

then we get

$$(f^* d\Gamma_G)_Q = (\delta(\phi|_{\mathcal{U}}))_Q = ((\delta\phi)|_{\mathcal{U}})_Q = (\delta\phi)_Q.$$

Hence, by recalling (5.1) and (5.3), it follows that

$$(\delta\phi)_Q = -(f^*(\Gamma_G \wedge \Gamma_G))_Q = -((f^*\Gamma_G) \wedge (f^*\Gamma_G))_Q = -(\phi \wedge \phi)_Q.$$

□

Theorem 5.2 and Proposition 3.3 yield immediately the following property.

**Corollary 5.1.** *Let $\mathcal{M}$ be a $M$-dimensional $C^2$ manifold and let $\phi \in Mat_L C^1 \mathcal{F}^1(\mathcal{M})$. Moreover, let $\mathcal{U} \subset \mathcal{M}$ be open and consider a $C^1$ map $f : \mathcal{U} \to G$. Then $(d\phi)_Q = -(\phi \wedge \phi)_Q$ for all $Q \in \mathcal{U} \cap \{f^*\Gamma_G = \phi|_\mathcal{U}\}^{(M+1)}$.*

Hence:

**Corollary 5.2.** *Let $\mathcal{M}$ be a $M$-dimensional $C^2$ manifold and let $\phi \in Mat_L C^1 \mathcal{F}^1(\mathcal{M})$ be such that $(d\phi)_P \neq -(\phi \wedge \phi)_P$ for a certain $P \in \mathcal{M}$. Then there exists a neighborhood $\mathcal{U}$ of $P$ such that $\mathcal{U} \cap \{f^*\Gamma_G = \phi|_\mathcal{U}\}^{(M+1)} = \emptyset$ for all $C^1$ maps $f : \mathcal{U} \to G$. In particular, if condition (5.2) is verified and $f : \mathcal{U} \to G$ is any $C^1$ map (with $\mathcal{U} \subset \mathcal{M}$ open), then one has $\mathcal{U} \cap \{f^*\Gamma_G = \phi|_\mathcal{U}\}^{(M+1)} = \emptyset$.*

University of Trento, Department of Mathematics, via Sommarive 14, 38123 Trento, Italy

*Email address*: `silvano.delladio@unitn.it`